\documentclass[10pt]{amsart}

\usepackage{eucal,url,amssymb,stmaryrd,enumerate,amscd,}
\usepackage[pagebackref,colorlinks=false]{hyperref}
\usepackage{amsfonts}
\usepackage{amsmath,amsthm,amssymb,amscd,enumerate,eucal,url,stmaryrd}
\usepackage{verbatim}

\numberwithin{equation}{section}
\newtheorem{thrm}{Theorem}[section]
\newtheorem{lemma}[thrm]{Lemma}
\newtheorem{prop}[thrm]{Proposition}
\newtheorem{cor}[thrm]{Corollary}

\newtheorem{rmrk}[thrm]{Remark}

\newtheorem{conv}[thrm]{Convention}
\begin{document}

\begin{abstract}
 We show that the fundamental 4-form on a quaternionic contact manifold of dimension at
 least eleven is closed { if and only if the torsion endomorphism} of
the Biquard connection  vanishes. This condition characterizes quaternionic contact
structures which are locally qc homothetic to 3-Sasakian structures.
\end{abstract}

\title[Quaternionic contact manifolds with {a closed fundamental} 4-form]%
{Quaternionic contact manifolds with {a closed fundamental} 4-form}
\date{\today}
\author{Stefan Ivanov}
\address[Stefan Ivanov]{University of Sofia, Faculty of Mathematics and
Informatics, blvd. James Bourchier 5, 1164, Sofia, Bulgaria}
\email{ivanovsp@fmi.uni-sofia.bg}

\author{Dimiter Vassilev}
\address[Dimiter Vassilev]{ Department of Mathematics and Statistics\\ University of New Mexico\\
Albuquerque, New Mexico, 87131}
\email{vassilev@math.unm.edu} \maketitle 

\setcounter{tocdepth}{2}

\section{Introduction}

A quaternionic contact (qc) structure, introduced in \cite{Biq1,Biq2},
appears naturally as the conformal boundary at infinity of the
quaternionic hyperbolic space. Such structures have been considered in connection with
the quaternionic contact Yamabe problem, \cite{Wei,IMV,IMV1}.
A particular case of this problem amounts to finding
the extremals and the best constant in the $L^2$ Folland-Stein
Sobolev-type embedding, \cite{F2} and \cite{FS}, on the quaternionic
Heisenberg group, see \cite{GV} and\cite{IMV1}.

{A  qc structure on a real (4n+3)-dimensional manifold $M$ is a codimension
three distribution $H$, { called the horizontal space}, locally given as the kernel of { a} 1-form $%
\eta=(\eta_1,\eta_2,\eta_3)$ with values in $\mathbb{R}^3$, {such
that,} the three 2-forms $d\eta_i|_H$ are the fundamental 2-forms
of a quaternionic structure on $H$. The 1-form $\eta$ is
determined up to a conformal factor and the action of $SO(3)$ on
$\mathbb{R}^3$. {Therefore} $H$ is equipped with a conformal class
$[g]$ of Riemannian metrics and a 2-sphere bundle of almost
complex structures, the quaternionic bundle $\bf{Q}$. The 2-sphere
bundle of one forms determines uniquely the associated metric and
a conformal change of the metric is equivalent to a conformal
change of the one forms. To every metric in the fixed conformal
class one can associate { a complementary to $H$ distribution $ V$
spanned by the Reeb vector fields $\xi_1,\xi_2,\xi_3$ and a linear
connection $\nabla$ preserving the qc structure and the splitting
$TM=H\oplus V$ provided $n>1$ \cite{Biq1}. This connection is
known as the Biquard connection.} The {qc Ricci tensor, the qc
scalar curvature $Scal$ of the Biquard connection  are obtained
from the curvature tensor by taking horizontal traces}.

The transformations preserving a given qc structure $\eta$, i.e.
$\bar\eta=\mu\Psi\cdot\eta$ for a positive smooth function $\mu$ and a
$SO(3)$ matrix $\Psi$ with smooth functions as entries, are called
\emph{quaternionic contact conformal (qc conformal)
transformations}. If the function $\mu$ is constant we have {\em
quaternionic contact homothetic (qc homothetic) transformations}.
The Biquard connection  is invariant under qc homothetic
transformations.

Examples of qc manifolds can be found in \cite{Biq1,Biq2,IMV,D1}.
In particular, any totally umbilic hypersurface of a quaternionic
K\"ahler or hyperk\"ahler manifold carries such a structure. An
extensiveley studied class of examples of quaternionic contact
structures are provided by the  3-Sasakian manifolds. The latter
can be defined as  $(4n+3)$-dimensional  (pseudo) Riemannian
manifold { of signature either (4n+3,0) or (4n,3)} whose
Riemannian cone is a hyperk\"ahler manifold { of signature
(4n+4,0) or (4n,4), respectively}. It was shown { in \cite{IMV}
that} the torsion endomorphism of the Biquard connection is the
obstruction for a given qc-structure to be locally qc homothethic
to a 3-Sasakian one provided the qc scalar curvature $Scal$ is not
identically zero. { Explicit examples of qc manifolds with zero or
non-zero torsion endomorphism were recently given in \cite{AFIV}}.
The quaternionic Heisenberg group, the quaternionic sphere of
dimension $4n+3$ with its standard 3-Sasakian structure and the qc
structures locally qc conformal to them  are characterized in
\cite{IV} by the vanishing of a tensor invariant, the qc-conformal
curvature defined in terms of the curvature and torsion of the
Biquard connection. Explicit examples
of non-qc conformally flat qc manifolds 
{are constructed} in \cite{AFIV}.

In this article we consider the  4-form $\Omega$ defining the $Sp(n)Sp(1)$
structure on the horizontal distribution and  { call it} {\it the fundamental
four-form.}

The purpose of the paper  is to show that
when the dimension of the manifold is greater than seven, the
fundamental { 4-form} form is closed if and only if the qc structure is
locally { qc homothetic to a 3-Sasakian one provided the qc scalar curvature does not vanish}.  We prove the following main result.

\begin{thrm}\label{main1}
Let $(M^{4n+3},\eta,\bf Q)$ be a $4n+3$-dimensional qc manifold.
For $n>1$ the following conditions are equivalent
\begin{enumerate}
\item[i)] The fundamental four form is closed, $d\Omega=0$;
\item[ii)] The torsion { endomorphism} of the Biquard
connection vanishes; \item[iii)] Each  Reeb vector field $\xi_l$,
defined in \eqref{bi1}, preserves the fundamental four form,
$\mathbb L_{\xi_l}\Omega=0$.
\end{enumerate}
Any of the above conditions imply that the qc scalar curvature is constant
and the vertical distribution is integrable.
\end{thrm}
Combining the last Theorem with  Theorem~1.3 and Theorem~7.11 in
\cite{IMV} {we obtain}
\begin{thrm}\label{main2}
Let $(M^{4n+3},\eta,\bf Q)$ be a $4n+3$-dimensional qc manifold.
For $n>1$ the following conditions are equivalent
\begin{enumerate}
\item[a)] {$(M^{4n+3},\eta,\bf{Q})$ has closed fundamental four
form, $d\Omega=0$}; \item[b)] { The torsion { endomorphism} of
the Biquard connection vanishes;} \item[c)] {$(M^{4n+3},g,\bf{Q})$
is a qc-Einstein manifold} (the trace-free part of the
qc Ricci tensor is zero); \item[d)] Each Reeb vector $\xi_l$
field preserves the horizontal metric and the quaternionic
structure simultaneously, $\mathbb L_{\xi_l}g=0,\quad \mathbb
L_{\xi_l}\bf Q\subset\bf Q$; \item[e)] Each  Reeb vector field
$\xi_l$ preserves the fundamental four form, $\mathbb
L_{\xi_l}\Omega=0$.
\end{enumerate}
\noindent If in addition the qc scalar curvature is { non-zero,
$Scal\not=0$}, then each of a), b), c), d) and e) is equivalent to
the following condition f).
\begin{enumerate}
\item[f)] $M^{4n+3}$ is locally qc homothetic to a 3-Sasakian
manifold, i.e., locally, there exists a $SO(3)$-matrix $\Psi$ with
smooth entries depending on an {auxiliary} parameter, such that,
the local { qc structure $(\epsilon\frac{Scal}{16n(n+2)}\Psi\cdot\eta,\bf
Q)$, $\epsilon=sign(Scal)$ is 3-Sasakian.}
\end{enumerate}
\end{thrm}

As an application of Theorem~\ref{main1} we give in the last section
 a proof of the equivalence of a) and f) in Theorem~\ref{main2}.
Thus, when the dimension of the qc manifold is greater than seven,
we establish  in a slightly different manner  Theorem~3.1 in \cite{IMV}.

\begin{rmrk}
On a seven dimensional  qc manifold, if the torsion endomorphism of the Biquard
connection vanishes then the fundamental four form is closed. We
do not know whether the converse holds or if there exists an
example of a seven dimensional qc manifold with a closed
fundamental four form and a non-vanishing torsion endomorphism.
This might be related to
 the well known fact that in dimension eight  an almost
quaternion hermitian structure with a closed fundamental  four
form is not necessarily quaternionic K\"ahler since Salamon
\cite{Sal} gave  a compact counter-example (see \cite{AFIV} for
non-compact counter-examples).
\end{rmrk}

\textbf{Organization of the paper.} The paper relies heavily on the
notion of Biquard connection introduced in \cite{Biq1} and the
properties of its torsion and curvature discovered in \cite{IMV}.
In order to make the present paper self-contained, in Section
\ref{s:review} we give a review of the notion of a quaternionic
contact structure and collect formulas and results from \cite{Biq1}
and \cite{IMV} that will be used.

\begin{conv}\label{conven}
\begin{enumerate}[a)]
\item We shall use $X,Y,Z,U$ to denote horizontal vector fields, i.e. $X,Y,Z,U\in H$;
\item $\{e_1,\dots,e_{4n}\}$ denotes an orthonormal basis of the horizontal
space $H$;
\item The triple $(i,j,k)$ denotes any cyclic permutation of
$(1,2,3)$.

\item $l$ and $m$ will be any numbers from the set $\{1,2,3\}, \quad
l,m \in\{1,2,3\}$.

\end{enumerate}
\end{conv}

\textbf{Acknowledgements} { We would like to thank the referee for
many valuable comments and remarks.} The research was done during
the visit of S.Ivanov to the Abdus Salam ICTP, Trieste as a Senior
Associate, Fall 2008. S.I. thanks  ICTP for providing the support
and an excellent research environment. S.I. is partially supported
by the Contract 082/2009 with the University of Sofia
`St.Kl.Ohridski'.
This work has been partially supported through  Contract ``Idei", DO 02-257/18.12.2008 and DID 02-39/21.12.2009

\section{Quaternionic contact manifolds}\label{s:qc str}

\label{s:review} In this section we will briefly review the basic notions of
quaternionic contact geometry and recall some results from \cite{Biq1} and
\cite{IMV}. For the purposes of this paper, a quaternionic contact (qc) manifold $(M, g,
\bf{Q})$ is a $4n+3$ dimensional manifold $M$ with a codimension three
distribution $H$ equipped with a metric $g$ and an Sp(n)Sp(1) structure,
i.e., we have
\begin{enumerate}
\item[i)] a 2-sphere bundle $\bf{Q}$ over $M$ of almost complex structures
$I_l\,:H \rightarrow H,\quad I_l^2\ =\ -1$, satisfying the commutation relations of the imaginary
quaternions $I_1I_2=-I_2I_1=I_3$ and
 ${\bf Q}= \{aI_1+bI_2+cI_3:\ a^2+b^2+c^2=1 \}$;
\item[ii)] $H$ is locally the kernel of a 1-form $\eta=(\eta_1,\eta_2,\eta_3)$ with
values in $\mathbb{R}^3$ satisfying the compatibility condition  \hspace{3mm}
$
2g(I_lX,Y)\ =\ d\eta_l(X,Y).
$
\end{enumerate}

The fundamental 2-forms $\omega_l$ of the quaternionic structure $\bf{Q}$ are {determined} by
\begin{equation}  \label{thirteen}
2\omega_{l|H}\ =\ \, d\eta_{l|H},\qquad \xi\lrcorner\omega_l=0,\quad \xi\in
V.
\end{equation}
The  4-form $\Omega$ defining the $Sp(n)Sp(1)$
structure on the horizontal distribution, called here the fundamental
four-form, is  defined (globally) on
the horizontal distribution $H$   by
\begin{equation}\label{fform}
\Omega=\omega_1\wedge\omega_1+\omega_2\wedge\omega_2+\omega_3\wedge\omega_3.
\end{equation}

On a quaternionic contact manifold there exists a canonical connection
defined in \cite{Biq1} when the dimension $(4n+3)>7$, and in \cite{D} in the
7-dimensional case.

\begin{thrm}
\cite{Biq1}\label{biqcon} {Let $(M, g,\bf{Q})$ be a quaternionic contact
manifold} of dimension $4n+3>7$ and a fixed metric $g$ on $H$ in the
conformal class $[g]$. Then there exists a unique connection $\nabla$ with
torsion $T$ on $M^{4n+3}$ and a unique supplementary subspace $V$ to $H$ in $%
TM$, such that:
\begin{enumerate}
\item[i)] $\nabla$ preserves the decomposition $H\oplus V$ and the $Sp(n)Sp(1)$-structure
on $H$;
\item[ii)] for $X,Y\in H$, one has $T(X,Y)=-[X,Y]_{|V}$;
\item[iii)] for $\xi\in V$, the endomorphism $T(\xi,.)_{|H}$ of $H$ lies in $%
(sp(n)\oplus sp(1))^{\bot}\subset gl(4n)$;
\end{enumerate}
\end{thrm}
We shall call the above connection \emph{the Biquard connection}.
Biquard \cite{Biq1} also described the supplementary subspace $V$.
{Locally, }$V$ is generated by vector fields $%
\{\xi_1,\xi_2,\xi_3\}$, such that
\begin{equation}  \label{bi1}
\begin{aligned} \eta_l(\xi_k)=\delta_{lk}, \qquad (\xi_l\lrcorner
d\eta_l)_{|H}=0,\quad (\xi_l\lrcorner d\eta_k)_{|H}=-(\xi_k\lrcorner
d\eta_l)_{|H}. \end{aligned}
\end{equation}
The vector fields $\xi_1,\xi_2,\xi_3$ are called Reeb vector fields or
fundamental vector fields.
{\ If the dimension of $M$ is seven, the conditions \eqref{bi1} do not
always hold. Duchemin shows in \cite{D} that if we assume, in addition, the
existence of Reeb vector fields as in \eqref{bi1}, then Theorem~\ref{biqcon}
holds. Henceforth, by a qc structure in dimension $7$ we shall 
mean a qc structure satisfying \eqref{bi1}.}

The torsion endomorphism $T_{\xi}=T(\xi,.) : H\rightarrow H, \quad \xi\in V,$
plays an important role in the qc geometry.
Decomposing the endomorphism $T_{\xi}\in(sp(n)+sp(1))^{\perp}$ into
symmetric part $T^0_{\xi}$ and skew-symmetric part $b_{\xi}$, $%
T_{\xi}=T^0_{\xi} + b_{\xi} $  Biquard shows in \cite{Biq1} that
$T_{\xi}$ is  completely trace-free, $tr\, T_{\xi}= tr\,
T_{\xi}\circ I=0,  I\in \bf{Q}$ and describes the properties of
the two components. Using the two $Sp(n)Sp(1)$-invariant
trace-free symmetric 2-tensors $T^0$ and  $U$ on $H$   {defined}
in \cite{IMV} by
\begin{gather*}
T^0(X,Y)\overset{def}{=}g((T_{\xi_1}^{0}I_1+T_{\xi_2}^{0}I_2+T_{%
\xi_3}^{0}I_3)X,Y),\quad U(X,Y)\overset{def}{=}-g(I_lb_{\xi_l}X,Y),
\end{gather*}
the properties of $T_{\xi}$ outlined in \cite{Biq1}  {give the
following identities, cf. \cite{IMV},}
\begin{gather}\label{propt}
T^0(X,Y)+T^0(I_1X,I_1Y)+T^0(I_2X,I_2Y)+T^0(I_3X,I_3Y)=0,
\\\label{propu} 3U(X,Y)-U(I_1X,I_1Y)-U(I_2X,I_2Y)-U(I_3X,I_3Y)=0.
\end{gather}
If $n=1$ then  $U$ vanishes identically, $U=0$, and the torsion is a
symmetric tensor, $T_{\xi}=T^0_{\xi}$.

The covariant derivatives with respect to the Biquard connection of
the almost complex structures and the vertical vectors are given by
\begin{equation}\label{der}
\nabla I_i=-\alpha_j\otimes I_k+\alpha_k\otimes I_j,\qquad
\nabla\xi_i=-\alpha_j\otimes\xi_k+\alpha_k\otimes\xi_j.
\end{equation}

{ It turns out that the vanishing of the $sp(1)$-connection
1-forms on $H$ is equivalent to the vanishing of the torsion
endomorphism of the Biquard connection, $T^0=U=0$ \cite{IMV}.

The first equation in \eqref{der} together with \eqref{fform}
imply that the fundamental four form is parallel with respect to
$\nabla$, $\nabla\Omega=0$ but it may not be closed because of the
torsion of the Biquard connection.}

Let $R=[\nabla,\nabla]-\nabla_{[\ ,\ ]}$ be the curvature tensor
of $\nabla$. We denote the curvature tensor of type (0,4) by the
same letter, $ R(A,B,C,D):=g(R(A,B)C,D),\  A,B,C,D \in
\Gamma(TM)$. The \emph{qc-Ricci forms} and
\emph{qc-scalar curvature} are defined  by ${4n\rho_l(A,B)=\sum_{a=1}^{4n}R(A,B,e_a,I_le_a)}, \quad
Scal=\sum_{a,b=1}^{4n}R(e_b,e_a,e_a,e_b),$ respectively. { It was shown in \cite{IMV} that
the $sp(1)$-part of $R$ is determined by the Ricci 2-forms and the
connection 1-forms by}
\begin{equation}  \label{sp1curv}
R(A,B,\xi_i,\xi_j)=2\rho_k(A,B)=(d\alpha_k+\alpha_i\wedge\alpha_j)(A,B),
\qquad A,B \in \Gamma(TM).
\end{equation}
{It is important to note that the horizontal part of the
Ricci 2-forms can be expressed in terms of the torsion of the
Biquard connection \cite{IMV}. For ease of reading, we collect
the necessary facts from \cite[Theorem~1.3, Theorem~3.12,
Corollary~3.14, Proposition~4.3 and Proposition~4.4]{IMV} with
slight modifications, using the equality
$4T^0(\xi_l,I_lX,Y)=T^0(X,Y)-T^0(I_lX,I_lY)$, and present them in
the form described in \cite{IV}.}
\begin{thrm}
\cite{IMV}\label{sixtyseven} 
On a $(4n+3)$-dimensional qc manifold, $n>1$  
the next formulas hold
\begin{equation*}
\begin{aligned} \rho_l(X,I_lY) \ & =\
-\frac12\Bigl[T^0(X,Y)+T^0(I_lX,I_lY)\Bigr]-2U(X,Y)-%
\frac{Scal}{8n(n+2)}g(X,Y),\\ T(\xi_{i},\xi_{j})& =-\frac
{Scal}{8n(n+2)}\xi_{k}-[\xi_{i},\xi_{j}]_{H}, \qquad Scal\  =\
-8n(n+2)g(T(\xi_1,\xi_2),\xi_3)
\\ T(\xi_i,\xi_j,X) &
=-\rho_k(I_iX,\xi_i)=-\rho_k(I_jX,\xi_j),\qquad
\rho_i(\xi_i,\xi_j)+\rho_k(\xi_k,\xi_j)=\frac1{16n(n+2)}\xi_j(Scal);\\
\rho_{i}(X,\xi_{i})\ & =\ -\frac {X(Scal)}{32n(n+2))} \ +\ \frac
12\, \left
(-\rho_{i}(\xi_{j},I_{k}X)+\rho_{j}(\xi_{k},I_{i}X)+\rho_{k}(\xi_{i},I_{j}X)
\right). \end{aligned}
\end{equation*}

In particular,  the vanishing of the horizontal trace-free part of
the Ricci forms is equivalent to the vanishing of the torsion
endomorphism of the Biquard connection. In this case the vertical
distribution is integrable, the qc scalar curvature is constant
and { if $Scal\not=0$ then}  the qc-structure is 3-Sasakian up to
a multiplication by a constant and an $SO(3)$-matrix with smooth
entries depending on an auxiliary parameter.
\end{thrm}
{ For the last part of the above Theorem we have adopted the
definition that a $4n+3$-dimensional (pseudo) Riemannian manifold
$(M,g_M)$ of signature either $(4n+3,0)$ or $(4n,3)$ has a
\emph{3-Sasakian structure} if the cone metric $t^2g_M+\epsilon dt^2$ on $M\times\mathbb
R$ is a hyperk\"ahler metric of signature  $(4n+4,0)$ (\textit{positive
3-Sasakian structure, $\epsilon=1$}) or (4n,4)  (\textit{negative
3-Sasakian structure, $\epsilon=-1$}), respectively, see \cite{Kon,Tano,Jel} for the negative case. In other words, the cone
metric has holonomy contained in $Sp(n+1)$ (see \cite{BGN}) or in
$Sp(n,1)$ (see \cite{AK1}), respectively.

We remind that, usually, a $4n+3$-dimensional Riemannian manifold
$(M,g)$ is called 3-Sasakian only in  the positive case, while the term pseudo 3-Sasakain is used in  the negative case. However,
we find it convenient to use the more general definition.

A 3-Sasakian manifold of  dimension $(4n+3)$ is Einstein
\cite{Kas,Tano,AK1}. If the metric is positive definite then it is
with positive Riemannian scalar curvature $(4n+2)(4n+3)$
\cite{Kas} and if complete it is compact with finite fundamental
group due to Myer's theorem. There are known many examples of
positive 3-Sasakian manifold of dimension $(4n+3)$, see \cite{BG}
and references therein  for a nice overview of positive 3-Sasakian
spaces. Certain $SO(3)$-bundles over quaternionic K\"ahler
manifolds with negative scalar curvature constructed in
\cite{Kon,Tano,Jel,AK1} supply examples of negative 3-Sasakian
manifolds. A natural definite metric on the negative 3-Sasakian manifolds is
considered in \cite{Tano,Jel} by changing the sign of the metric
on the vertical $SO(3)$-factor. With respect to this metric the
negative 3-Sasakian manifold becomes an $A$-manifold in the
terminology of \cite{Gray}, its Riemannian Ricci tensor has
precisely two constant eigenvalues, $-4n-14$ (of multiplicity
$4n$) and  $4n+2$ (of multiplicity $3$) see \cite{Jel}, and the
Riemannian scalar curvature is the negative constant
$-16n^2-44n+6$ \cite{Tano,Jel}. Explicit examples of negative
3-Sasakian manifolds are constructed in \cite{AFIV}.

\section{Local structure equations of qc manifolds}\label{s:local str eqs}

We derive the local structure equations of a qc structure in
terms of the $sp(1)$-connection forms of the Biquard connection and
the qc scalar curvature.

\begin{prop}\label{str}
Let $(M^{4n+3},\eta,\bf Q)$ be a (4n+3)- dimensional qc manifold
with qc scalar curvature $Scal$. Let $s=\frac{Scal}{8n(n+2)}$ be the
normalized qc scalar curvature. The following  equations hold
\begin{gather}\label{streq}
2\omega_i=d\eta_i+\eta_j\wedge\alpha_k-\eta_k\wedge\alpha_j +s
\eta_j\wedge\eta_k,\\\label{str2}
d\omega_i=\omega_j\wedge(\alpha_k+s\eta_k)-\omega_k\wedge(\alpha_j+s\eta_j)-\rho_k\wedge\eta_j+
\rho_j\wedge\eta_k+\frac12ds\wedge\eta_j\wedge\eta_k,\\\label{strom}
d\Omega=\sum_{(ijk)}\Big[2\eta_i\wedge(\rho_k\wedge\omega_j-\rho_j\wedge\omega_k)+
ds\wedge\omega_i\wedge\eta_j\wedge\eta_k\Big],
\end{gather}
where $\alpha_l$ are the $sp(1)$-connection 1-forms of the Biquard
connection, $\rho_l$ are the Ricci 2-forms and $\sum_{(ijk)}$ is the
cyclic sum of even permutations of $\{1,2,3\}$.

In particular,   for a 3-Sasakian manifold the structure equations
have the form
\begin{equation}\label{33sas}
d\eta_i=2\omega_i+2\epsilon\eta_j\wedge\eta_k
\end{equation}
and the normalized qc scalar curvature is $s=2\epsilon$, {
where $\epsilon= 1$  if the 3-Sasakian structure is positive and
$\epsilon=-1$ in the negative 3-Sasakian case.}
\end{prop}
\begin{proof}
From the definition \eqref{thirteen} of the fundamental 2-forms
$\omega_l$ we have (see also \cite{IMV})
\begin{equation}\label{fifteen}
2\omega_m=(d\eta_m)_{|_H}=d\eta_m -
\sum_{l=1}^3{\eta_l\wedge({\xi_l}\lrcorner d\eta_{m})} + \sum_{1\leq
l<p\leq 3}{d\eta_m(\xi_l,\xi_p) \eta_l \wedge \eta_p}.
\end{equation}
It is shown in \cite{Biq1} that the $sp(1)$-connection 1-forms
$\alpha_l$ on $H$ are given by
\begin{gather}  \label{coneforms}
\alpha_i(X)=d\eta_k(\xi_j,X)=-d\eta_j(\xi_k,X), \quad X\in H, \quad
\xi_i\in V.
\end{gather}
The $sp(1)$-connection 1-forms $\alpha_l$ on the vertical space
$V$ were determined in \cite{IMV}:
\begin{multline}  \label{coneform1}
\alpha_i(\xi_l)\ =\ d\eta_l(\xi_j,\xi_k) -\
\delta_{il}\left(\frac{Scal}{ 16n(n+2)}\ +\ \frac12\,\left(\,
d\eta_1(\xi_2,\xi_3)\ +\ d\eta_2(\xi_3,\xi_1)\ + \
d\eta_3(\xi_1,\xi_2)\right)\right).
\end{multline}
A straightforward calculation
using \eqref{coneforms} and \eqref{coneform1} gives the equivalence of \eqref{fifteen} and \eqref{streq}.
Taking the exterior derivative of \eqref{streq}, followed by an application of
\eqref{streq} and \eqref{sp1curv} implies \eqref{str2}. The last formula, \eqref{strom},
follows from \eqref{str2} and  definition \eqref{fform}.

 { For the last part of the theorem consider the cone $N=M^{4n+3}\times\mathbb R^+$ equipped with
 a} (pseudo) { almost hyperhermitian} structure $(G_N,\phi_l)$ where
$G_N$ is a (pseudo) Riemannian metric of signature (4n+4,0) (resp.
(4n,4)) for $\epsilon=1$ (resp. $\epsilon=-1$) and $\phi_l$,
$l=1,2,3$, are three anti-commuting almost complex structures. The
1-form $dt$ on $\mathbb R^+$ and the three almost complex
structures are related to three 1-forms $\eta_l$ on $M^{4n+3}$
defined by $\eta_l=\epsilon\,\phi_l\, (\frac1t dt)$, where we used
the same notation for  both a tensor and its lift to a  tensor on
the tangent bundle of $N$ identifying $M$ with the slice $t=1$ of
$N$. We may write the metric $G_N$ and the three K\"ahler 2-forms
on $N$ as follows:
\begin{equation}\label{e:cone metric and 2-forms}
\begin{aligned}
 G_N=t^2g+\epsilon t^2 ( \eta_1^2 + \eta_2^2 + \eta_3^2)  +\epsilon dt^2;\qquad
 F_i^\epsilon= t^2 \omega_i +\epsilon  t^2\eta_j\wedge \eta_k - t\eta_i\wedge dt,
\end{aligned}
\end{equation}
 where $g={G_N}_{|_H}$ and $H=\cap_{l=1}^3 Ker\, \eta_l$. A qc structure on $M^{4n+3}$ is defined
by the three 1-forms $\eta_l$ \cite{Biq1}.

It is straightforward to check from the second equation in
\eqref{e:cone metric and 2-forms} that the 2-forms $F_l^\epsilon$
are closed precisely when \eqref{33sas} holds. Therefore the cone
metric $G_N$  is hyperk\"ahler, i.e. $M^{4n+3}$ is 3-Sasakian, if
and only if \eqref{33sas} is fulfilled  due to Hitchin's theorem
\cite{Hit}, { which is valid with the same proof in the case of
non-positive definite metrics}.

To compute the qc scalar curvature of a 3-Sasakian manifold, we
use equations \eqref{33sas} to find $ {{\xi_i}\lrcorner
d\eta_{j}}_{|H}=0,\ d\eta_i(\xi_j,\xi_k)=2\epsilon,\
d\eta_i(\xi_i,\xi_k)=d\eta_i(\xi_i,\xi_j)=0. $   We calculate from
\eqref{coneforms} and \eqref{coneform1} the $sp(1)$-connection
1-forms
$
\alpha_l=-\Big(\frac{Scal}{16n(n+2)}+\epsilon\Big)\eta_l. $ The
last identity and \eqref{sp1curv} yield
$
\rho_l(X,Y)=\frac12d\alpha_l(X,Y)=-\Big(\frac{Scal}{16n(n+2)}+\epsilon\Big)\omega_l(X,Y),
$ which compared with the first equation in
Theorem~\ref{sixtyseven} gives  $T=U=0, \ Scal=16n(n+2)\epsilon$,
see \cite{IMV}.

\end{proof}

\section{Proof of Theorem~\ref{main1}}
First we show that if $T^0=U=0$, then $d\Omega=0$. Indeed, in this
case, Theorem~\ref{sixtyseven} implies
\begin{equation}\label{e:rhota for Omega closed}
\rho_l(X,Y)=-s\omega_l(X,Y), \quad \rho_l(\xi_m,X)=0, \quad
\rho_i(\xi_i,\xi_j)+\rho_k(\xi_k,\xi_j)=0,
\end{equation}
since $Scal$ is constant and the horizontal distribution is
integrable. Using the just obtained identities in \eqref{strom}
gives $d\Omega=0$ { which proves the implication
$ii)\rightarrow i)$}.

{ To finish the proof of the theorem we shall apply the next Lemma.}
\begin{lemma}\label{lemdom}
On a qc manifold of dimension $(4n+3)>7$ we have the identities
\begin{gather}\label{domu}
U(X,Y)=-\frac1{16n}\sum_{a=1}^{4n}\Big[d\Omega(\xi_i,X,I_kY,e_a,I_je_a)+d\Omega(\xi_i,I_iX,I_jY,e_a,I_je_a)
\Big]\\\label{domt}
T^0(X,Y)=\frac1{8(1-n)}\sum_{(ijk)}\sum_{a=1}^{4n}\Big[d\Omega(\xi_i,X,I_kY,e_a,I_je_a)-d\Omega(\xi_i,I_iX,I_jY,e_a,I_je_a)
\Big].
\end{gather}
\end{lemma}
\begin{proof}
Equation \eqref{strom} together with the  first equality in
Theorem~\ref{sixtyseven} yield
\begin{equation}\label{rom}
d\Omega(\xi_i,X,I_kY,e_a,I_je_a)=4(n-1)\rho^0_k(X,I_kY)+
2\rho^0_j(X,I_jY)-2\rho^0_j(I_iX,I_kY),
\end{equation}
where $\rho^0$ is the horizontal trace-free part of $\rho$ given by
\begin{equation}\label{trr}
\rho^0_l(X,I_lY) \  =\
-\frac12\Bigl[T^0(X,Y)+T^0(I_lX,I_lY)\Bigr]-2U(X,Y).
\end{equation}
A substitution of \eqref{trr} in \eqref{rom}, combined with  the properties of
the torsion, \eqref{propt} and \eqref{propu} give
\begin{equation}\label{lasp}
2(n-1)\Big[T^0(X,Y)+T^0(I_kX,I_kY)\Big]
+8nU(X,Y)=-\sum_{a=1}^{4n}d\Omega(\xi_i,X,I_kY,e_a,I_je_a).
\end{equation}
{Applying again \eqref{propt} and \eqref{propu} in \eqref{lasp} we
see that} $U$ and $T^0$ satisfy \eqref{domu} and \eqref{domt},
respectively which proves the lemma.
\end{proof}

The well known Cartan formula yields $\mathbb
L_{\xi_l}\Omega=\xi_l\lrcorner
d\Omega+d(\xi_l\lrcorner\Omega)=\xi_l\lrcorner d\Omega, $ since
$\Omega$ is horizontal. The latter formula together with the
already proved implication $ii)\rightarrow i)$ and
Lemma~\ref{lemdom} complete the proof of Theorem~\ref{main1}.

From Lemma~\ref{lemdom} we easily derive
\begin{cor}\label{lemdomu}
If  one of the Reeb vector fields preserves the fundamental four
form on a qc manifold of dimension $(4n+3)>7$ then $U=0$ and the
torsion endomorphism of the Biquard connection is symmetric, $T_{\xi_l}=T_{\xi_l}^0$.
\end{cor}

\section{Proof  of  Theorem~\ref{main2}}

 In this section we give the proof of the equivalence of parts a)
 and f) of Theorem~\ref{main2}. The remaining claims follow from Theorem \ref{main1} and
\cite[Theorem~1.3 and Theorem~7.11]{IMV}. The idea is the same as
in the proof of Theorem~3.1 in \cite{IMV}, namely we show that
both a) and f) are equivalent to the fact that the cone over $M$
is locally hyperk\"ahler. However, here, the proof is based on the
fundamental 4-form. In one direction, let $d\Omega=0$.
Theorem~\ref{main1} implies that the torsion of the Biquard
connection vanishes, while Theorem~\ref{sixtyseven} {shows} that
the qc scalar curvature is constant and the vertical distribution
is integrable. Let $Scal\not=0$ and $\epsilon=sign(Scal)$. The qc
structure $\eta'=\epsilon\frac {Scal}{16n(n+2)}\eta$ has
normalized qc scalar curvature $s'=2\epsilon$ and $d\Omega'=0$.
For simplicity, we shall denote $\eta'$ with $\eta$ and, in fact,
drop the $'$ everywhere.

In the first step of the proof we show that the cone
$N=M\times\mathbb R^+$ with the {structure $(G_N,F_l^{\epsilon})$ defined by \eqref{e:cone metric and 2-forms}
 has holonomy contained in $Sp(n+1)$ {for
$\epsilon>0$ and in $Sp(n,1)$ for $\epsilon<0$}. To this end we
consider the following four form on $N$
\begin{equation}\label{con4}
    F=F^{\epsilon}_1\wedge F^{\epsilon}_1+F^{\epsilon}_2\wedge F^{\epsilon}_2+F^{\epsilon}_3\wedge F^{\epsilon}_3.
\end{equation}
Applying \eqref{streq}, \eqref{str2} and
\eqref{strom}, we calculate
\begin{multline}\label{nnnn}
dF_i=tdt\wedge(2\omega_i+2\epsilon\eta_j\wedge\eta_k-
d\eta_i)+t^2d(\omega_i+\epsilon\eta_j\wedge\eta_k)\\
=t\, dt\wedge (2\omega_i+2\epsilon\eta_j\wedge\eta_k+\eta_j\wedge\alpha_k-\eta_k\wedge\alpha_k+s\eta_j\wedge\eta_k)
+\epsilon t^2 (2\omega_j+\eta_i\wedge\alpha_k)\wedge\eta_k-\epsilon t^2 (2\omega_k-\eta_i\wedge\alpha_j)\wedge\eta_j)\\
+t^2\left (\omega_j\wedge(\alpha_k+s\eta_k)-\omega_k\wedge(\alpha_j+s\eta_j)-\rho_k\wedge\eta_j+\rho_j\wedge\eta_k+\frac 12 ds\wedge\eta_j\wedge\eta_k \right).
\end{multline}
 A short computation, using \eqref{streq}, \eqref{str2},
\eqref{strom} and \eqref{nnnn}, gives}
\begin{multline}
dF=2\sum_{i=1}^3 dF_i\wedge F_i=t^4\sum_{(ijk)} \left [2\eta_i\wedge(\omega_j\wedge\rho_k - \omega_k\wedge \rho_j) +ds\wedge\omega_i\wedge\eta_j\wedge\eta_k\right ]\\
+t^3dt\wedge\sum_{(ijk)} \left [4\epsilon\omega_i\wedge\eta_k\wedge\eta_j-2s\omega_i\wedge\eta_j\wedge\eta_k - 4\rho_k\wedge\eta_i\wedge\eta_j
 -ds\wedge\eta_i\wedge\eta_j\wedge\eta_k\right ]\\
 =2t^4d\Omega - 4t^3\sum_{(ijk)}dt\wedge(\rho_i+2\epsilon\omega_i)\wedge\eta_j\wedge\eta_k=0,
\end{multline}
taking into account the first equality in
Theorem~\ref{sixtyseven}, \eqref{e:rhota for Omega closed} and
$s=2\epsilon$, which hold when $d\Omega=0$ by Theorem \ref{main1}.

Hence, $dF=0$ and the holonomy of the cone metric is contained
either in $Sp(n+1)Sp(1)$ or in $Sp(n,1)Sp(1)$ provided $n>1$
\cite{AS},  i.e. the cone is quaternionic K\"ahler manifold
provided $n>1$. Note that when $n=1$ this conclusion can not be
reached  in the positive definite case due to the 8-dimensional
compact  counter-example constructed by S. Salamon \cite{Sal}
(for non compact counter-examples see \cite{AFIV}).

It is a classical result (see e.g \cite{Bes} and references
therein) that a quaternionic K\"ahler manifolds of dimension
bigger than four (of arbitrary signature) are Einstein. This fact
implies that the cone $N=M\times \mathbb R^+$ with the  metric
$g_N$ must be Ricci flat (see e.g. \cite[p.267]{Bes}) and
therefore it is locally hyperk\"ahler { since the $sp(1)$-part
of the Riemannian curvature vanishes and therefore it can be
trivialized locally by a parallel sections} (see e.g.
\cite[p.397]{Bes}). This means that locally there exists a
$SO(3)$-matrix $\Psi$ with smooth entries, possibly depending on
$t$,  such that the triple of two forms $(\tilde F_1,\tilde
F_2,\tilde F_3)= \Psi\cdot (F_1,F_2,F_3)^T$ consists of closed
2-forms constitute the fundanental 2-forms of the local
hyperk\"ahler structure. Consequently $(M,\Psi\cdot\eta)$ is
locally a 3-Sasakian manifold.

The fact that f) implies a) is trivial since the 4-form $\Omega$
is invariant under rotations and rescales by a constant when the
metric on the horizontal space $H$ is replaced by a homothetic to
it metric.

\begin{rmrk} It follows from the above discussion that the cone over a (4n+3)-dimensional qc
manifold carries either a $Sp(n+1)Sp(1)$ or a $Sp(n,1)Sp(1)$
structure which is closed exactly when $d\Omega=0$ provided $n>1$.
\end{rmrk}

\begin{rmrk} {
An example of a qc structure satisfying $T^0=U=Scal=0$ can be obtained as follows.
Let $M^{4n}$ be a hyperk\"ahler
manifold with closed and locally exact K\"ahler forms
$\omega_l=d\eta_l$. The total space of an ${\mathbb R^3}$-bundle over the hyperk\"ahler
manifold $M^{4n}$ with connection 1-forms $\eta_l$ is an example of a qc structure with $T^0=U=Scal=0$.} The qc
structure is determined by the three 1-forms $\eta_l$ satisfying
$d\eta_l=\omega_l$ which yield $T^0=U=Scal=0$. In particular, the
quaternionic Heisenberg group which locally is the unique qc
structure with flat Biquard connection on $H$, see \cite{IV}, can
be considered as an $\mathbb R^3$ bundle over a 4n-dimensional
flat hyperk\"ahler $\mathbb R^{4n}$. A compact example is provided
by a $T^3$-bundle over a compact hyperk\"ahler manifold $M^{4n}$
such that each closed K\"ahler form $\omega_l$ represents integral
cohomology classes. Indeed, since $[\omega_l]$, $1\leq l\leq 3$
define integral cohomology classes on $M^{4n}$, the well-known
result of Kobayashi \cite{Kob} implies that there exists a circle
bundle $S^1 \hookrightarrow M^{4n+1} \to M^{4n}$, with connection
$1$-form $\eta_1$ on $M^{4n+1}$ whose curvature form is $d\eta_1 =
\omega_1$. Because $\omega_l$ $(l=2,3)$ defines an integral
cohomology class on $M^{4n+1}$, there exists a principal circle
bundle $S^1 \hookrightarrow M^{4n+2} \to M^{4n+1}$ corresponding
to [$\omega_2$] and a connection $1$-form $\eta_2$ on $M^{4n+2}$
such that $\omega_2=d\eta_2$ is the curvature form of $\eta_2$.
Using again the result of Kobayashi, one gets a $T^3$-bundle over
$M^{4n}$ whose total space has a qc structure  satisfying
$d\eta_l=\omega_l$ which yield $T^0=U=Scal=0$.

We do not know
whether there are other examples satisfying the
conditions $T^0=U=Scal=0$.
\end{rmrk}

\end{document}